\documentclass[11pt]{article}
\usepackage{graphics,color,amsmath,amsthm,amsfonts,verbatim}

 \DeclareMathOperator{\lca}{lca}
 \DeclareMathOperator{\diam}{diam}

\newcommand{\M}{\mathcal{M}}

\newcommand{\A}{\mathcal{A}}
\renewcommand{\L}{\mathcal{L}}
\newcommand{\T}{\mathcal{T}}
\newcommand{\R}{\mathbb{R}}
\newcommand{\e}{\varepsilon}

\theoremstyle{plain}
 \newtheorem{theorem}{Theorem}
 \newtheorem{lemma}{Lemma}

 \newtheorem{proposition}[lemma]{Proposition}

 \theoremstyle{definition}

\begin{document}

\title{Limitations to Fr\'echet's Metric Embedding Method}
\date{}

\author{Yair Bartal\thanks{Supported in part by a grant from the
Israeli National Science Foundation} \and Nathan
Linial\thanks{Supported in part by a grant from the Israeli
National Science Foundation} \and Manor Mendel\thanks{Supported in
part by the Landau Center.} \and Assaf Naor }

\maketitle

\begin{abstract}
Fr\'echet's classical isometric embedding argument has evolved to
become a major tool in the study of metric spaces. An important
example of a Fr\'echet embedding is Bourgain's embedding
\cite{bourgainembedding}. The authors have recently shown
\cite{blmn1} that for every $\e>0$ any $n$-point metric space
contains a subset of size at least $n^{1-\e}$ which embeds into
$\ell_2$ with distortion $O\bigl(\frac{\log(2/\e)}{\e}\bigr)$. The
embedding used in \cite{blmn1} is non-Fr\'echet, and the purpose
of this note is to show that this is not coincidental.
Specifically, for every $\epsilon >0$,
we construct arbitrarily large $n$-point metric spaces,
such that the distortion of any Fr\'echet embedding into $\ell_p$ on subsets of
size at least $n^{1/2 + \epsilon}$ is
$\Omega\left((\log n)^{1/p}\right)$.
\end{abstract}

\section{Introduction}

\newcommand{\skel}[1]{\mathrm{skel}({#1})}
\newcommand{\leaves}[1]{\mathrm{leaves}({#1})}

Given two metric spaces $(X,d_X)$, $(Y,d_Y)$ and an embedding
$f:X\to Y$, the {\em distortion} of $f$ is defined as:
\[
\mathrm{dist}(f)=\sup_{\substack{x,y\in X\\x\neq
y}}\frac{d_Y(f(x),f(y))}{d_X(x,y)}\cdot \sup_{\substack{x,y\in X\\x\neq
y}}\frac{d_X(x,y)}{d_Y(f(x),f(y))}.
\]
We denote by $c_Y(X)$ the least distortion with which $X$ may be
embedded in $Y$. When $c_Y(X)\le \alpha$ we say that $X$
$\alpha$-embeds into $Y$ and denote
$X\stackrel{\alpha}{\hookrightarrow}Y$. When there is a bijection
$f$ between two metric spaces $X$ and $Y$ with
$\mathrm{dist(f)}\le \alpha$ we say that $X$ and $Y$ are
$\alpha$-equivalent. For a class of metric spaces $\M$, $c_\M(X)$
is the minimum $\alpha$ such that $X$ $\alpha$-embeds into some
metric space in $\M$. For $p\ge 1$ we denote $c_{\ell_p}(X)$ by
$c_p(X)$. The parameter $c_2(X)$ is known as the {\em Euclidean
distortion} of $X$. A fundamental result of Bourgain
\cite{bourgainembedding} states that $c_2(X)=O(\log n)$ for every
$n$-point metric space $(X,d)$.

For a general metric space $(X,d)$ with no additional a-priori
structure, there is a dearth of genuinely "interesting"
constructions of Lipschitz mappings on $X$. One significant
exception to this rule is provided by the distance functions
$x\mapsto d(x,A)$ for some $\emptyset\neq A\subset X$. Of course,
we can generate more examples by constructing Lipschitz functions
to any finite dimensional normed space, the coordinates of which
are distance functions. Observe that a mapping $f:X\to
\ell_{\infty}$ is $L$-Lipschitz if and only if each of its
coordinates is $L$-Lipschitz. These facts were put to good use in
the classical observation of Fr\'echet that every metric space
isometrically embeds into some $\ell_{\infty}(\Gamma)$ (see
\cite{benlin, matbook}). Fr\'echet's embedding only uses distance
functions for singleton $A$'s, but more sophisticated refinements
of this basic idea have appeared over the years. In what follows
we call an embedding $f:X\to \R^{2^V\setminus\{\emptyset\}}$ a
{\em Fr\'echet embedding} if for every $A\subset X$ there is
$\alpha_A\in \R$ such that $f(x)=(\alpha_A d(x, A))_{A\in
2^X\setminus\{\emptyset\}}$.

Bourgain's embedding of finite metric spaces in
$\ell_2$~\cite{bourgainembedding} is an instance of a Fr\'echet
embedding in which the coefficients $\alpha_A \in [0,1]$ depends
only on the cardinality of $A$. Bourgain's probabilistic method of
producing a good Fr\'echet embedding has subsequently found many
applications \cite{jls,llr,matexpander,matdist,feige,tz}.

As stated in the abstract, the present note is motivated
by our recent Ramsey-type result:

\begin{theorem}[\cite{blmn1}]\label{thm:lowerlarge}
There exists an absolute constant $C>0$ such that for every
$\alpha>2$, any finite metric $M$ contains a subset $N\subset M$
for which $|N|\geq |M|^{1-C\frac{\log\alpha}{\alpha}}$ and
$c_2(N)\leq \alpha$.
\end{theorem}

The embedding used in Theorem~\ref{thm:lowerlarge} is not a
Fr\'echet embedding. In view of the past success of Fr\'echet
embeddings, and in particular Bourgain's embedding which gives the
asymptotically best possible bound for embedding the whole metric
space, it is natural to ask whether this (by now standard) method
is applicable to Ramsey-type problems. This note provides a
negative answer to this question. We find a certain range of the
parameters, for which Fr\'echet embeddings fail to achieve tight
bounds for Ramsey type questions:

\begin{theorem} \label{thm:bourgain fails} For every $1/2<\delta\le 1$ there is a constant $C(\delta)>0$ such that
for infinitely many integers $n$ there are $n$-point metric spaces $M_n$,
such that for any Fr\'echet embedding $f:M_n\rightarrow \ell_p$, and any
subset $V$ of $M_n$ of size at least $n^{\delta}$, $\text{dist}(f|_V)\geq
(C(\delta)\cdot\log n)^{1/p}$.
\end{theorem}

We end the paper with a short discussion in which we comment on
the embedding used in Theorem~\ref{thm:lowerlarge}, showing that
it can be viewed as a different natural generalization of
Bourgain's embedding.

\section{The Construction}

Theorem~\ref{thm:bourgain fails} is proved by exhibiting an
explicit example of an unbounded family of metric spaces for which
every Fr\'echet embedding fails to yield the appropriate Ramsey
type result. The example, denoted $Z=Z_{k,\e,h,m}$, is an
amalgamation of two types of metric spaces:
\begin{itemize}
\item
The elements of the first metric space, called $X=\L_{h,k}$, are  the leaves of
a complete binary tree $\T_h$ of height $h$. The metric on the leaves is
defined by
 \[ d_X(x_1,c_2)= k^{-l(\lca(x_1,x_2))} .\]
Here $\lca(x_1,x_2)$ denotes the least common ancestor of $x_1$
and $x_2$ in $\T_h$ and $l(u)$ is the depth of the vertex $u$ in
$T$. The parameter $k> 1$ will be be fixed later. A crucial
property of X is that for any $x,y,z\in X$, if
$d_X(x,y)<d_X(x,z)$, then
\begin{equation}
 d_X(x,z)=d_X(y,z)\geq k d_X(x,y)
\end{equation}
(Such a metric space is called a $k$-HST \cite{bartal1}).

\item The second metric space, $Y=Y_{m,\e}$, is the one-dimensional
metric on the points
 \[ Y=  \left\{y_0=\e, y_1=\e\left(1+\tfrac{1}{4}\right), y_2=\e\left(1+\tfrac{1}{4}+\tfrac{1}{16}\right),
   \ldots, y_{m-1}=\e\left( \sum_{i=0}^{m-1} 4^{-i}  \right) \right\} .\]
\end{itemize}

The points set of $Z$ is $X\times Y$,  hence its size is $2^h m$.
The distance is defined by
\[
d_Z((x,y),(x',y'))= \begin{cases} d_Y(y,y') & x=x' \\
  d_Y(y,y_0)+d_X(x,x')+d_Y(y_0,y') & x\neq x' .\end{cases}
\]

\begin{figure}[t]
\input{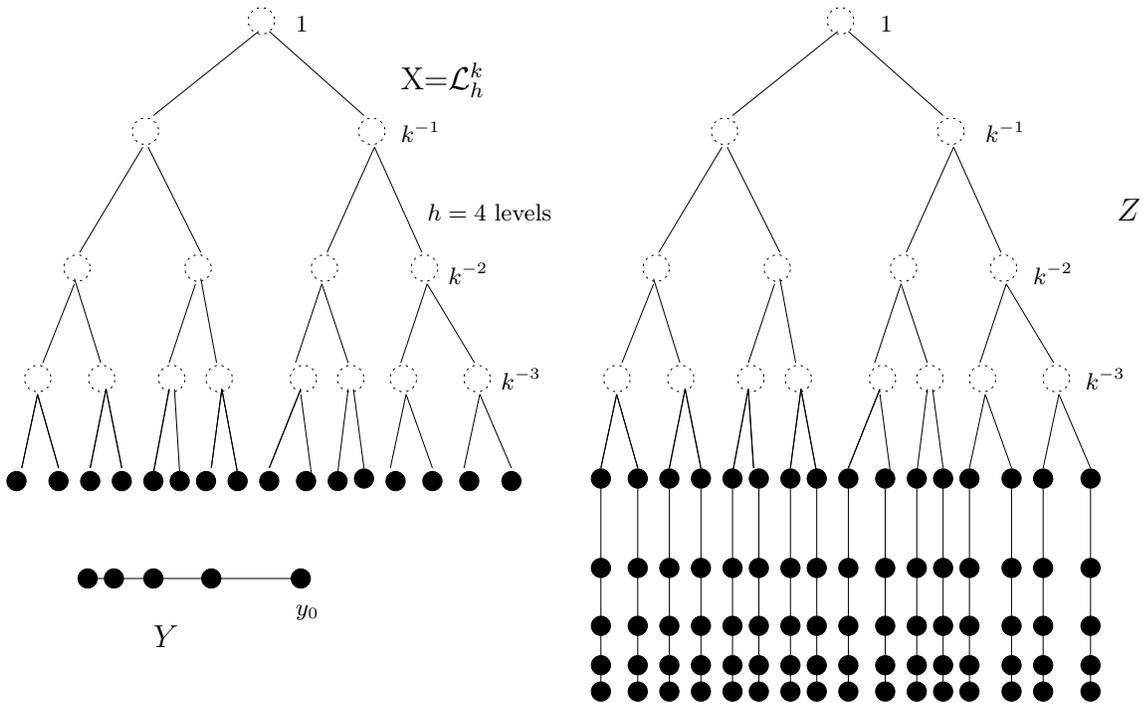}
\caption{The metric space $Z$ serving as
the example in Theorem~\ref{thm:bourgain fails}.}
\label{fig:example}
\end{figure}

A schematic description of $Z$ is given in
Figure~\ref{fig:example}.

Here is a sketch of our proof. A preliminary
step (Lemma~\ref{lem:sub-subset}) shows that all sufficiently
large (of size $\ge n^{\frac{1}{2}+\epsilon}$) subsets $V \subseteq Z$
have a structure that is similar to the whole of $Z$. Namely,
there is a large subset $U \subseteq V$ that spans a complete
binary subtree (of a slight modification) of $\cal T$. Moreover,
in each copy of $Y$ that $U$ meets, it contains at least two elements.

From this point on, we may assume, then, that we are dealing with
a space that is similar to the original $Z$ with two
modifications: The depth has shrunk to one half its original
value, (since we moved from $n$ to $n^{\frac{1}{2}+\epsilon}$
points) and attached to each leaf is a {\em tail} i.e., a
two-point metric space. On the other hand, we are now considering
the whole of this (sub)space, and not subsets thereof. (We
actually cannot ignore the points in $Z\setminus U$, which may
participate in the subsets to which Fr\'echet embedding assigns
nonzero weights. However, this is only a technicality and can be
ignored at this level of discussion.)

Recall that we are dealing with Fr\'echet embeddings $f$ that are
defined by assigning weights $\alpha_A$ to subsets $A$. Consider a
tail $\{y_1,y_2\}$, and suppose that $A$ is disjoint from the copy
of $Y$ to which $y_1,y_2$ belong. Then $|d(y_1,A) - d(y_2,A)| =
d(y_1,y_2)$, entailing a lower bound on $\mathrm{Lip}(f)$. In
other words, if large weights $\alpha_A$ are allotted to subsets
$A$ that miss many tails, then $\mathrm{Lip}(f)$ is big. This is
made precise in Lemma~\ref{lem:lower beta}, which gives a lower
bound on $\mathrm{Lip}(f)$ in terms of the numbers $\alpha_A$ and
$\zeta_A$ - the fraction of tails that $A$ misses.

It follows that in order to keep $\mathrm{Lip}(f)$ small, we ought to
put significant weights $\alpha_A$ on subsets
with large $\zeta_A$ i.e., those that meet many tails.
This, however, tends to increase $\mathrm{Lip}(f^{-1})$, as we now explain.
Let $z_1$ and $z_2$ belong to two distinct tails. In order for
$|d(z_1,A) - d(z_2,A)|$ to be nearly equal to $d(z_1,z_2)$, the
set $A$ must meet exactly one of the two tails containing the $z_i$
(in which case we say that $A$ {\em splits} the least common
ancestor of $z_1,z_2$). It is not hard to see
(Equation~\ref{eq:induction tree}) that if $\zeta_A$ is large,
then $A$ can split only a few vertices, and these necessarily reside
far from the root. Note also that if $A$ fails to split
$\lca(z_1,z_2)$, then
$|d(z_1,A) - d(z_2,A)| \le \frac{1}{k} d(z_1,z_2)$,
which is significant, since $k$ is large.
The precise argument is made in Lemma~\ref{lem:lower inverse}.

Omitting some additional technicalities (Lemma~\ref{lem:lower gamma}), these two considerations
can be traded off against each other
to yield the desired result.

\subsubsection*{Some notations and definitions concerning trees:}

\begin{enumerate}
 \item Let ${\cal T}$ be a tree. We denote its root by $r({\cal T})$
and its set of leaves by $\mathrm{leaves}({\cal T})$.
For a subset $A\subset \mathrm{leaves}({\cal T})$ let ${\cal
T}(A)$ be the subtree spanned by $A$, i.e. the union of all the simple paths
between elements of $A$ and $\lca(A)$.
The subtree $\mathcal{T}_u$ rooted at $u$ consists of the
union of all the monotone paths in ${\cal T}$ between $u$ and its descendants.
We recall that $l_{\cal T}(u)$ is the depth of $u$, i.e. its distance from
$r({\cal T})$.

\item Let $\mathcal{T}=(V,E)$ be a tree and suppose that $\mathcal{T}$ contains a vertex $u\neq r(\mathcal{T})$
with exactly two neighbors $v$ and $w$. We {\it simplify}
$\mathcal{T}$ by removing the vertex $u$ and adding the
edge $(v,w)$. Another type of simplification step takes place,
if $r(\mathcal{T})$ has a single neighbor $v$. In this case,
we remove $r(\mathcal{T})$ and declare $v$ as the root.
The skeleton skel$(\mathcal{T})$ is the
(uniquely defined) tree obtained from $\mathcal{T}$ by
carrying out all possible simplifications.

\item For a subset $U\subset Z$, we denote $U_X=\{x\in X;\; (\{x\}\times Y)
\cap U \neq \emptyset \}$.
\end{enumerate}

\begin{lemma}\label{lem:sub-subset} Let $h,m\ge 2$ be integers, and $Z=Z_{k,\e,h,m}$.
Let $V\subset Z$ be such that $|V|>2m+2^h$. Denote
\[ \rho= \left( \frac{|V|-2^h}{m} \right )^{1/h} -1.\]
Then there is a subset $U\subset V$ such that for every $x\in
U_X$, $|(\{x\} \times Y) \cap U|\ge 2$ and $\skel{{\cal
T}_h(U_X)}$ is a complete binary tree of depth at least $
\frac{\rho}{6\log (2/\rho)} h$.
\end{lemma}
\begin{proof}
Denote $S=\{x\in X;\ |V\cap (\{x\}\times Y )|\ge 2\}$. Since $|S|m+2^h\ge
|V|$, we have that
$$
|S|\ge \frac{|V| -2^h}{m}=(1+\rho)^h > 2.
$$

Set $0<\eta=\frac{\rho}{2+\rho}<0.5$. We will first construct
inductively a subtree ${\cal T'}$ of ${\cal T}$, as follows. Let
$r$ be the root of ${\cal T}$ and $u_1,u_2$ its two children.
Recursively, assume that we have constructed ${\cal T'}_{u_1}$,
${\cal T'}_{u_2}$. ${\cal T'}$ is constructed according to the
following rule: if $\min\{|\leaves{{\cal T}_{u_1}}\cap
S|,|\leaves{{\cal T}_{u_2}}\cap S| \}\ge \eta |S|$ then ${\cal
T'}$ is obtained by attaching ${\cal T'}_{u_1}$ and ${\cal
T'}_{u_2}$ to $r$. Otherwise, assume without loss of generality
that $|\leaves{{\cal T}_{u_1}}\cap S|> (1-\eta)|S|$. In this case
${\cal T'}$ is obtained by attaching only ${\cal T}_{u_1}$ to $r$.
Note that in the previous step of this procedure, ${\cal
T'}_{u_i}$ were constructed with respect to sets $\leaves{{\cal
T}_{u_i}}\cap S$.

Let $u$ be a leaf of $\skel{\cal T'}$ and denote its distance from
the root of $\skel{\cal T'}$ by $\lambda$. Consider the path in ${\cal
T}$ between $r$ (the root of ${\cal T}$) and $u$. Clearly $\lambda$ is
precisely the number of times along this path that the first
option was taken in the above construction. We prove by induction on $h$, the height
of $\cal T$,  that
$|S|\eta^\lambda (1-\eta)^{h-\lambda}\le 1$.
The case $h=0$ is trivial, so assume
that $h>0$. Let $r,u_1,u_2$ be as above, and assume that $u_1$
belongs to the path connecting $r$ and $u$. In the first case of
the above construction, we get by the induction hypothesis that:
\[
1\ge |\leaves{{\cal T}_{u_1}}\cap
S|\eta^{\lambda -1}(1-\eta)^{h-1-(\lambda -1)}\ge |S|\eta^\lambda(1-\eta)^{h- \lambda}.
\]

In the second case of the above construction we get
\[
1\ge |\leaves{{\cal T}_{u_1}}\cap
S|\eta^{\lambda }(1-\eta)^{h-1-\lambda }\ge |S|\eta^\lambda(1-\eta)^{h- \lambda}.
\]

Let $\hat{h}$ be the minimal distance between a leaf of
$\skel{\cal T'}$ and its root. From the discussion above
(and noting that $\rho<1$) it follows that:
$$
\hat{h}\ge \frac{\log|S|+h\log
(1-\eta)}{\log\left(\frac{1-\eta}{\eta}\right)}\ge \frac{ h \log
(1+\rho) + h \log (\frac{2}{2+\rho})}{\log \left(\frac{2}{\rho}
\right )} \ge \frac{\rho}{6\log(2/\rho)} h.
$$

Consider all the vertices of depth exactly $\hat{h}$ in
$\skel{\cal T'}$. The tree $\skel{\mathcal{T}'}$ truncated at
these points forms a complete binary tree of height $\hat{h}$. For
each such vertex $u$ we choose some $s_u\in S\cap {\cal T'}_u$.
The set $S'=\{s_u;\text{ $u$ has depth }\hat{h}\text{ in
}\skel{\cal T'}\}$ is easily seen to have the required properties.
\end{proof}

Our next two lemmas provide a lower bound on the expansion of any Fr\'echet embedding $Z \hookrightarrow \ell_p$,
using (almost) only the structure of $Y$.

Let $f: Z\to \ell_p^{\,2^{Z}\setminus\{\emptyset\}}$, $1\leq p<\infty$,
be a Fr\'echet embedding given by $f(x)_A=\alpha_A d(x,A)$,
let $V \subseteq Z$ be a large set as in Theorem~\ref{thm:bourgain fails}
and Let $U\subset V$ be a subset as in Lemma~\ref{lem:sub-subset}.
For every $\emptyset \neq A\subset Z$ write:
$$
\zeta_A=\frac{|\{x\in U_X;\ (\{x\}\times Y)\cap A= \emptyset\}|}{|U_X|}.
$$
Namely, $\zeta_A$ is the fraction of of the tails in $U$ that $A$ misses.
Define
\begin{equation}\label{eq:beta}
\beta(f,U)=\left(\sum_{\emptyset \neq A\subset Z}|\alpha_A|^p
\zeta_A \right )^{1/p}.
\end{equation}
Informally, $\beta(f)$ is a weighted average of the fraction of
tails $\{x\}\times Y$,  $x\in U$, that are missed by a subset of
$Z$. If a subset $A$ is disjoint from a tail $\{x\}\times Y$, that
contains the points $(x,u),(x,v)$, then
$|d_Z((x,u),A)-d_Z((x,v),A)| = d_Y(u,v)$. Consequently, we can
expect to obtain lower bound on the Lipschitz constant of $f$ in
terms of $\beta$. Indeed:

\begin{lemma}\label{lem:lower beta}
Let $U\subset Z $ be a subset as in  Lemma~\ref{lem:sub-subset}.
Then $\mathrm{Lip}(f|_U)\ge \beta(f,U)$.
\end{lemma}
\begin{proof} We may assume without loss of generality that $\alpha_A\ge 0$ for every $\emptyset\neq A\subset Z$.
Set $\beta=\beta(f,U)$. Note that
\[
\frac{1}{|U_X|}\sum_{x\in U_X}\left(\sum_{A:\; A\cap
Z_x=\emptyset} \alpha_A^p\right) =\sum_{\emptyset \neq A\subset
Z}\alpha_A^p\frac{|\{x\in U_X;\ A\cap (\{x\} \times Y)
=\emptyset\}|}{|U_X|}=\beta^p.
\]
Hence there is some $x\in U_X$ such that $\sum_{A:\; A\cap Z_x=\emptyset}
\alpha_A^p\ge \beta^p$. As mentioned, if
$\emptyset\neq A\subset Z$ is disjoint from $\{x\}\times Y$, then
$d_Z((x,y),A)=d_Y(y,y_0)+d_Z((x,y_0),A)$ for every $(x,y) \in Z$.
Recall that for $x\in U_X$, there are $y\neq y'$ such that $(x,y),(x,y') \in (\{x\}\times Y) \cap U$. Therefore
\begin{align*}
\|f(x,y)-f(x,y')\|_p^p &=\sum_{A\subset Z} \alpha_A^p
|d_Z((x,y),A)-
d_Z((x,y'),A)|^p\\
&\geq \sum_{\substack{A\subset Z\\
A \cap (\{x\}\times Y)=\emptyset}} \alpha_A^p |d_Y(y_0,y)-d_Y(y_0,y')|^p\geq
\beta^p d_Z((x,y),(x,y'))^p,
\end{align*}
so that $\mathrm{Lip}(f|_U)\ge \beta$.
\end{proof}

Similar to the definition of $\beta(f,U)$, we also define
\begin{equation} \label{eq:gamma}
\gamma(f,U)=\left (\sum_{\emptyset\neq A\not \supset
U}\alpha_A^p\right )^{1/p}.
\end{equation}
In words, $\gamma(f,U)$ is the weighted average of weights
assigned to subsets that do not contain $U$. The metric of $Y$,
was purposely chosen so that no two distances are equal.
Consequently, a subset $A\not\supset U$ makes a nonzero
contribution to $\mathrm{Lip}(f|_U)$. Indeed:

\begin{lemma}\label{lem:lower gamma}
Let $U \subseteq Z$ be as in Lemma~\ref{lem:lower beta}. Assume
that $k\geq 4$, and $k^{-h}\geq \e$. Then $ \mathrm{Lip}(f|_U)\ge
4^{-m} 2^{-\frac{|X|\cdot |Y|}{p}}\gamma(f,U)$.
\end{lemma}
\begin{proof}
We may assume again that $\alpha_A\ge 0$ for every
$\emptyset \neq A\subset Z$. Let $A_0$ be such that
$\alpha_{A_0}=\max\{\alpha_A;\ \emptyset\neq A\not \supset U\}$.
Clearly:
\[
\alpha_{A_0}^p\ge 2^{-|X|\cdot |Y|}\sum_{\emptyset\neq A\not
\supset V}\alpha_A^p.
\]
Since $A_0\not \supset U$, there is some $(x,u)\in U\setminus A_0$. By
assumption, $U$ contains some $(x,v)$ with $v\neq u$. It is
enough to show that $|d_Z((x,u),A_0)-d_Z((x,v),A_0)|\ge
4^{-m}\cdot\mathrm{diam}(Y)$, since then:
$$
\|f(x,u)-f(x,v)\|_p\ge
\alpha_{A_0}|d_Z((x,u),A_0)-d_Z((x,v),A_0)|\ge
4^{-m}\alpha_{A_0}\mathrm{diam}(Y)\ge 4^{-m}\alpha_{A_0}d_Y(u,v),
$$
and the result will follow from the lower estimate for
$\alpha_{A_0}$.

To verify that $|d_Z((x,u),A_0)-d_Z((x,v),A_0)|\ge
4^{-m}\cdot\mathrm{diam}(Y)$, distinguish between the possible
points where the distances $d_Z((x,u),A_0)$ and $d_Z((x,v),A_0)$
are attained. Let $(x',a)\in A_0$, $(x'',b)\in A_0$ satisfy
$d_Z((x,u), A_0)=d_Z((x,u),(x',a))$, $d_Z((x,v),
A_0)=d_Z((x,v),(x'',b))$. Observe that since
$k^{-h}>\e/3>\diam(Y)$, if $x\in\{x',x''\}$ then $x'=x''=x$, in
which case $d_Z((x,u), A_0)=d_Y(u,a)$, $d_Z((x,v), A_0)=d_Y(v,b)$,
and we conclude since the definition of the metric on $Y$ implies
that $|d_Y(u,a)-d_Y(u,b)|\ge 4^{-m+1}\e>4^{-m}\diam(Y)$.
Therefore, assuming that $x\notin \{x',x''\}$, $d_Z((x,u),
A_0)=d_Y(u,y_0)+d_X(x,x')+d_Y(y_0,a)$ and $d_Z((x,v),
A_0)=d_Y(v,y_0)+d_X(x,x'')+d_Y(y_0,b)$. If $d_X(x,x')=d_X(x,x'')$,
we use the fact that if $a,b,u,v \in Y$ but $\{a,b\}\neq \{u,v\}$
then $\left|[d_Y(u,y_0)+d_Y(y_0,a)]-[d_Y(v,y_0)+d_Y(y_0,b)]\right|
\geq 4^{-m}\diam(Y)$. On the other hand, if $d_X(x,x')\neq
d_X(x,x'')$ then $|d_X(x,x') - d_X(x,x'')| \geq k^{-h} \geq \e$,
and so
\begin{eqnarray*}
|d_Z((x,u),A_0)-d_Z((x,v),A_0)|&\ge&
|d_X(x,x')-d_X(x,x'')|-\\&\phantom{\le}&|d_Y(u,y_0)-d_Y(v,y_0)|-|d_Y(a,y_0)-d_Y(b,y_0)|\\
&\ge& \e-2\mathrm{diam}(Y)\ge \mathrm{diam}(Y).
\end{eqnarray*}
\end{proof}

Our next lemma uses the special structure of $X={\cal L}_{h,k}$ to
bound from below the inverse of the contraction,
$\mathrm{Lip}((f|_U)^{-1})$, in terms of $\beta(f,U)$ and
$\gamma(f,U)$. In what follows we always use the convention that
$\mathrm{Lip}(f^{-1})=\infty$ if $f$ is non injective.

\begin{lemma}\label{lem:lower inverse}
Let $U\subset Z$ be as in Lemma~\ref{lem:sub-subset} .
Then for every Fr\'echet
embedding $f:Z\to \ell_p^{2^Z\setminus \{\emptyset\}}$,
\[
\mathrm{Lip}\left((f|_U)^{-1}\right)\ge
\left[\frac{2\beta(f,U)^p}{\hat{h}}+\left(\frac{2+2k^h\diam(Y)}{k}\right)^p\gamma(f,U)^p\right]^{-1/p}.
\]
\end{lemma}
\begin{proof} We use a shorthand notation ${\cal T}={\cal T}_h$,
$\gamma=\gamma(f,U)$ and $\beta=\beta(f,U)$. We assume as usual
that $\alpha_A\ge 0$ for every $\emptyset \neq A\subset Z$. Denote
by ${\cal T'}={\cal T}_h(U_X)$ the subtree generated by $U_X$. By
our assumption $\skel{{\cal T}'}$ is a complete binary tree of
depth $\hat {h}$. Let $u; u_1, u_2$ be a vertex in $\skel{{\cal
T}'}$ and its two children in ${\cal T}$. We say that a subset
$A\subset Z$ \emph{splits} $u$ if $A$ intersects exactly one of
the sets $\mathrm{leaves}({\cal T}_{u_1})\times Y
,\mathrm{leaves}({\cal T}_{u_2})\times Y$.

For $u\in\skel{\cal T'}$, denote by $l(u)=l_{\skel{\cal T'}}(u)$
the depth of $u$ in $\skel{{\cal T'}}$.
We first claim that for $\emptyset\neq A\subset Z$,
\begin{equation}\label{eq:induction tree}
\sum_{\substack{u\in \skel{\cal T'}\\
A \text{ splits }u}}2^{-    l(u)}\leq 2\zeta_A.
\end{equation}
The proof is by induction on $\hat{h}$, the height of $\skel{\cal
T'}$. Let $r$ be the root of $\skel{\cal T'}$. Denote by $u_1,u_2$
the children of $r$ in ${\cal T}$, and by $v_1$, $v_2$ the
children of $r$ in $\skel{\cal T'}$. Set $\zeta=\zeta_A$ and
define $\zeta_1,\zeta_2$ by:
$$
\zeta_i=\frac{|\{x\in \leaves{{\cal T'}_{v_i}}; (\{x\}\times Y)
\cap A= \emptyset\}|}{|\leaves{{\cal T'}_{v_i}}|},\quad i=1,2.
$$
Note that since $\skel{\cal T'}$ is a complete binary tree
$\zeta=(\zeta_1+\zeta_2)/2$. If $A$ does not split $r$ then by
induction:
\begin{align*}
\sum_{\substack{u\in \skel{\cal T'}\\ A \text{ splits }u}} 2^{-
l(u)}&= \frac{1}{2}\left ( \sum_{\substack{u\in \skel{{\cal
T'}_{v_1}}\\ A \text{ splits }u}} 2^{-( l(u)-1)}
+\sum_{\substack{u\in \skel{{\cal T'}_{v_2}}\\ A
\text{ splits }u}} 2^{-(    l(u)-1)} \right )\\
 &\leq \frac{1}{2} (2\zeta_1 + 2\zeta_2)=2\zeta .
\end{align*}

On the other hand, if $A$ splits $r$, it does not intersect one of
$\leaves{{\cal T}_{u_1}}\times Y$, $\leaves{{\cal T}_{u_2}}\times
Y $, say it does not intersect
$\leaves{{\cal T}_{u_2}}\times Y$. In this case $\zeta_2=1$ so
that by the induction hypothesis:
\begin{align*}
\sum_{\substack{u\in \skel{\cal T'}\\ A \text{ splits }u}}
2^{-    l(u)}=
 1+ \frac{1}{2} \left(
 \sum_{\substack{u\in \skel{{\cal T'}_{v_1}}\\ A \text{ splits }u}} 2^{-(   l(u)-1)}
  \right ) \leq 1 + \frac{1}{2} 2\zeta_1 = \zeta_2+\zeta_1 = 2\zeta.
\end{align*}
This finishes the proof of (\ref{eq:induction tree}). Now, by the
definition of $\beta=\beta(f,U)$,

\begin{align*}
\sum_{l=0}^{\hat{h}-1} \sum_{\substack{u\in \skel{\cal T'}\\ u \text{ in depth }l\\
\text{in } \skel{\cal T'}}} 2^{-l} \sum_{\substack{\emptyset\neq
A\subset Z\\A \text{ splits }u}} \alpha_A^p = \sum_{\emptyset\neq
A\subset Z} \alpha_A^p \sum_{\substack{u\in \skel{\cal T'}\\ A
\text{ splits }u}}2^{-  l(u)} \leq \sum_{\emptyset\neq A\subset Z}
\alpha_A^p 2\zeta_A = 2 \beta^p.
\end{align*}
It follows that there exists some $l\in\{0,\ldots,\hat{h}-1\}$
such that
\[ 2^{-l} \sum_{\substack{u\in \skel{\cal T'}\\ u
       \text{ in depth }l\\\text{in }\skel{\cal T'}}}
 \left ( \sum_{\substack{\emptyset\neq A\subset Z\\A \text{ splits }u}} \alpha_A^p\right )  \leq
\frac{2\beta^p}{\hat{h}}.\]
So there exists a vertex $u\in
\skel{\cal T'}$ (in depth $l$ in $\skel{\cal T'}$)  such that
\[  \sum_{\substack{\emptyset\neq A\subset Z\\A \text{ splits }u}} \alpha_A^p  \leq
\frac{2\beta^p}{\hat{h}}.\]

Denote by $u_1,u_2$ the two children of $u$ {\bf in} ${\cal T}$. Since $u$ is
a vertex of $\skel{\cal T'}$ there are $x_1\in \leaves{{\cal T}_{u_1}}\cap
U_X$ and $x_2\in \leaves{{\cal T}_{u_2}}\cap U_X$ (in particular
$\lca(x_1,x_2)=u$). Fix $(x_1,y_1)\in  U$ and $(x_2,y_2)\in  U$. The
observation here is that if a subset $A\subset Z$ does not split $u$ then one
of the following two cases happens:
\begin{enumerate}
\item $A \cap ({\cal T}_u \times Y) = \emptyset$. In this case, $d_Z((x_1,y_1),A)= d_Z((x_2,y_2),A)$, and therefore
\[ |d_Z((x_1,y_1),A)- d_Z((x_2,y_2),A)|=0.\]
\item
For $i=1,2$
there is $x_i'\in \leaves{{\cal T}_{u_i}}$ and $y_i',y_i''\in Y$
such that
$d_Z((x_i,y_i),A)=d_Y(y_i,y_i')+d_X(x_i,x_i')+d_Y(y_0,y_i'')$.
Hence,
{\setlength\arraycolsep{1pt}
\begin{eqnarray*}
|d_Z((x_1,y_1),A)- d_Z((x_2,y_2),A)|&\leq&
d_X(x_1,x_1')+d_X(x_2,x_2')+2\diam(Y)\\
&\le& \frac{2}{k}d_X(x_1,x_2)+2\diam(Y)\\
&\le& \frac{2+2k^h\diam(Y)}{k}d_Z((x_1,y_1),(x_2,y_2)).
\end{eqnarray*}}
\end{enumerate}

Therefore {\setlength\arraycolsep{1pt}
\begin{eqnarray*}
 \|f((x_1,y_1))&-&f((x_2,y_2))\|^p_p =\\
 &=& \sum_{\substack{\emptyset\neq A\subset Z\\A \text{ splits }u}} \alpha_A^p
 |d_Z( (x_1,y_1),A)-d_Z((x_2,y_2),A)|^p\\
 &\phantom{\le}&+\sum_{\substack{\emptyset\neq A\subset Z\\A \text{ doesn't splits }u}}
 \alpha_A^p|d_Z( (x_1,y_1),A)-d_Z((x_2,y_2),A)|^p\\
 &\le&[d_Z((x_1,y_1),(x_2,y_2))]^p\sum_{\substack{\emptyset\neq A\subset Z\\A
 \text{ splits }u}} \alpha_A^p\\
 &\phantom{\le}&+\left(\frac{2+2k^h\diam(Y)}{k}d_Z((x_1,y_1),(x_2,y_2))\right)^p\sum_{\emptyset\neq A\not \supset
U}\alpha_A^p\\
&\le&
\left[\frac{2\beta^p}{\hat{h}}+\left(\frac{2+2k^h\diam(Y)}{k}\right)^p\gamma^p\right]
[d_Z((x_1,y_1),(x_2,y_2))]^p.
\end{eqnarray*}}
\end{proof}

\medskip

\begin{proof}[Proof of Theorem \ref{thm:bourgain fails}.] Let $h,m\ge 2$
where $\epsilon = k^{-h}/2$
Define $n=m2^h$, so that $|Z|=n$.


Let $f:Z\to \ell_p^{2^Z\setminus \{\emptyset\}}$ be a Fr\'echet
embedding. Fix $0.5<\delta<1$. We can always choose $h$ and $m$
such that $n^\delta=\Omega(2^{h}+mn^{2\delta-1})$ ($2^h$ of order
$n^\delta$ and $m$ of order $n^{1-\delta}$). Fix such $h$ and $m$,
and let $V\subset Z$ be such that $|V|\ge n^\delta$. By
Lemma~\ref{lem:sub-subset} there is a subset $U\subset V$ which
satisfies the conditions of Lemma \ref{lem:lower inverse} with
$\hat{h}\geq C(\delta)\log n$. Set $\beta=\beta(f,U)$,
$\gamma=\gamma(f,U)$. It follows from Lemmas \ref{lem:lower beta},
\ref{lem:lower gamma}, \ref{lem:lower inverse} that:
$$
\mathrm{dist}(f|_V)\ge \mathrm{dist}(f|_U)\ge c\frac{\max\{\beta,
2^{-n/p}4^{-m}\gamma\}}{\max\left\{\frac{\beta}{(C(\delta)\log
n)^{1/p}},\frac{\gamma}{k}\right\}},
$$
for some universal constant $c>0$.

If $\frac{\beta}{(C(\delta)\log n)^{1/p}}\ge\frac{\gamma}{k}$ then
we deduce that $\mathrm{dist}(f|_V)\ge (C(\delta)\log n)^{1/p}$,
as required. Otherwise we get the lower bound
$\mathrm{dist}(f|_V)\ge c k2^{-n/p}4^{-m}$. Recall that we are
still free to choose $k$, so that the required result also follows
from this case provided that $k$ is large enough.
\end{proof}

\section{Discussion}

The goal of Theorem~\ref{thm:bourgain fails} has been to provide
an example for which Fr\'echet embeddings fail to achieve the best
possible bounds. This is done for the problem of embedding subsets
of size at least $n^\delta$ into $\ell_p$, where $\delta > 1/2$.
We suspect that the same holds also for $\delta \le 1/2$. However,
the example presented here does not immediately apply to this
case.

Finally we comment on a concept due to Matou\v{s}ek and
Rabinovich~\cite{yuri-jiri} which is a different natural
generalization of Bourgain's embedding. Given a finite metric
space $(V,\rho)$, we say that a one-dimensional metric $\sigma$ on
$V$ is {\em dominated} by $\rho$ if $\rho(x,y) \ge \sigma(x,y)$
for every $x,y \in V$. The polytope $\ell_1^{\text{dom}}(\rho)$ is
the convex hull of all one-dimensional metrics dominated by
$\rho$. It is natural to ask whether embeddings into
$\ell_1^{\text{dom}}$ can be used for Ramsey-type problems. We
observe that this is indeed the case. The following theorem is a
consequence of the main Ramsey-type theorem of \cite{blmn1}:


\begin{theorem}
\label{l1dom_ramsey}
For every finite metric $M$, and every
$\alpha>2$, there exists a subset $N\subset M$
of cardinality $\geq |M|^{1-C\frac{\log\alpha}{\alpha}}$
that is $O(\alpha)$-equivalent to some metric in $\ell_1^{\mathrm{dom}}(N)$.
Here $C>0$ is an absolute constant.
\end{theorem}

The proof of Theorem~\ref{l1dom_ramsey} is in two steps: (i) We
recall that Theorem~\ref{thm:lowerlarge} follows from a
Ramsey-type theorem of \cite{blmn1} where the target metric space
for the embedding is an {\em ultrametric} (which is isometrically
embeddable in $\ell_2$). (ii) Every ultrametric $\rho$ is
$O(1)$-equivalent to some metric in $\ell_1^{\mathrm{dom}}(\rho)$.
This is a simple fact, a proof of which, is sketched below.

\begin{proposition}
Every ultrametric $\rho$ is $16$-equivalent to some metric in $\ell_1^{\mathrm{dom}}(\rho)$.
\end{proposition}
\begin{proof}[Sketch of a proof.]
As shown in \cite{blmn1}, $\rho$ has a $4$-embedding in some
$4$-HST. So, it suffices to prove that $\rho$ is 4 equivalent to
some metric in $\ell_1^{\mathrm{dom}}(\rho)$ for the case where
$\rho$ itself is a 4-HST metric. Let us recall: A $k$-HST is
defined by a tree $T$ in which every internal vertex $v$ is
assigned a weight $\Delta(v) > 0$. If $v$ is a child of $u$ then
$\Delta(u) \ge k \Delta(v)$. The metric is defined of $T$'s leaves
via $d(x,y)=\Delta(\lca(x,y))$. The metric in
$\ell_1^{\mathrm{dom}}(\rho)$ that is $4$-equivalent to $\rho$ is
constructed through a probabilistic argument. Associate with every
edge $e=uv$ a weight $\epsilon_e \in \{-\frac{3}{8},
\frac{3}{8}\}$. This is done uniformly and independently over all
edges. Associated with this is a mapping
$\varphi=\varphi_{\epsilon}$ of $T$'s leaves to the real line,
$\varphi(x) = \sum \Delta(u) \epsilon_{uv}$. This sum extends over
all edges $uv$ in the directed path from $T$'s root to $x$. It is
easy to verify that the one dimensional metric induced on the
leaves is dominated by $\rho$. This is now averaged over all
possible choices of the $\epsilon_e$'s.
\end{proof}

\bibliographystyle{plain}
\bibliography{frechet}

\bigskip
\bigskip

\noindent Yair Bartal, Institute of Computer Science, Hebrew
University, Jerusalem 91904, Israel. \\{\bf yair@cs.huji.ac.il}

\medskip
\noindent Nathan Linial, Institute of Computer Science, Hebrew
University, Jerusalem 91904, Israel. \\ {\bf nati@cs.huji.ac.il}

\medskip
\noindent Manor Mendel, Institute of Computer Science, Hebrew
University, Jerusalem 91904, Israel. \\ {\bf
mendelma@cs.huji.ac.il}

\medskip
\noindent Assaf Naor, Theory Group, Microsoft Research, One
Microsoft Way 113/2131, Redmond WA 98052-6399, USA. \\ {\bf
anaor@microsoft.com}

\bigskip
\bigskip

\end{document}